\documentclass[12pt]{amsart}
\usepackage{euscript}
\usepackage{multicol}
\usepackage{amssymb}
\usepackage{amsmath}
\usepackage{times} 

\newtheorem{theorem}{Theorem}[section]
\newtheorem{proposition}{Proposition}[section]

\newtheorem{lemma}{Lemma}[section]

\numberwithin{equation}{section}
\begin{document}
\title[]
{ Lower  bounds on the eigenvalue sums of the Schr\"odinger operator and the spectral conservation law}
\author[O. Safronov]{ Oleg Safronov}
\address{Department  of  Mathematics and Statistics, UNCC, Charlotte, NC 28223, USA}
\email{osafrono@uncc.edu}
\thanks{The  author would  like to thank Andrej Zlatos, Rupert Frank  and Alexander Gordon
for useful discussions and remarks}

\maketitle


\section{Statement of the main results}

In this paper we  consider the  Schr\"odinger operator
$$
H=-\Delta-V(x),\quad V>0,
$$
acting in the  space $L^2({\Bbb R}^d)$.  We  study the relation between the behavior  of $V$ at the infinity and the properties of the negative spectrum of $H$. According to the  Cwikel-Lieb-Rozenblum estimate \cite{CLR}, \cite{19},\cite{22} the number of negative eigenvalues
of $H$ satisfies the relation
\begin{equation}\label{clr}
N\leq C \int V^{d/2}dx,\qquad d\geq 3,
\end{equation}
where $C$ is independent of $V$. Similar estimates hold in dimensions $d=1$ and $d=2$.
Therefore, if $V$ decays at the infinity fast enough, then $N$ is finite.
The question arise if finiteness of $N$ implies a qualified  decay of $V>0$ at the infinity?

 One can try to formulate  theorems whose assumptions
 contain as less as possible information about $V$. But then it is not clear how to define $H$.
To keep our arguments simple,  we shall assume that
$$
V\in L^\infty({\Bbb R}^d).
$$
 The  first  result which is related to the two dimensional case is proven in \cite{GNL}.
\begin{theorem}\label{N} {\bf (Grigoryan-Netrusov-Yau \cite{GNL} )} Let $d=2$ and
let $V>0$ be a bounded function on ${\Bbb R}^2$. Assume  that the negative spectrum of $H$ consists of  $N$  eigenvalues. Then
$V\in L^1({\Bbb R}^2)$ and
\begin{equation}\label{V<N1}\int_{{\Bbb R}^2} Vdx\leq c_0 N.
\end{equation}
where  the constant $c_0$ does not depend on $V$.
\end{theorem}
  Instead of proving \eqref{V<N1}, we will establish the estimate
\begin{equation}\label{V<N2}\int_{{\Bbb R}^2} Vdx\leq (6^4+12+4N)N.
\end{equation}
While the relation \eqref{V<N2} does  not give any new  information compared to \eqref{V<N1}, we think that 
 comprehension of  the arguments  related to this inequality will  help to understand the idea of Theorem~\ref{main}.

One is tempted  to say that  a straitforward   generalization of this result to the case $d\geq3$ should establish finiteness of the integral  $\int V^{d/2}dx$ under the condition that $N<\infty$. But such a   statement is false (see Theorem~\ref{smolchan}), because 
$$
\frac{(d-2)^2}4\int_{{\Bbb R}^d}\frac {|u|^2}{ |x|^2}dx\leq \int_{{\Bbb R}^d} |\nabla u|^2 dx,\qquad \forall u\in C^\infty_0({\Bbb R}^d \setminus \{0\}),\,\,\, d\geq3,
$$ 
which  means that  the operator $-\Delta-V$ with the potential  $V(x)=4^{-1}(d-2)^2/(1+|x|^2)$
does not have negative eigenvalues.
Instead, 
one can easily prove the implication
\begin{equation*}
N<\infty
\implies
\int_{|x|>1} \frac{V(x)}{|x|^{d-2+\varepsilon}}dx<\infty,\qquad \forall \varepsilon>0.
\end{equation*}
Indeed, without loss of generality, we can assume that $H\geq0$ on the complement  of a large ball $B$. It  means that 
$$
\int V(x)|\phi(x)|^2dx\leq \int |\nabla \phi(x)|^2dx,\qquad \forall \phi\in C_0^\infty({\Bbb R}^d\setminus B).
$$
It remains to take $\phi(x)= |x|^{-d/2+1-\varepsilon/2}$ for $|x|$ large enough.

Theorem~\ref{N}   deals with the  number of negative eigenvalues $\lambda_n$, whereas  one can study similar questions
related to the  sums
$$
\sum_n |\lambda_n|^p,\qquad p>0.
$$
\begin{theorem}\label{LTd3} {\bf (Damanik-Remling \cite{DR})} Let $B_1$ be the unit  
ball in  ${\Bbb R}^d$ and 
let $V>0$ be a bounded function on ${\Bbb R}^d \setminus B_1$.
 Assume  that the negative spectrum of the operator $H=-\Delta-V$ with  the Neumann  boundary condition  on the unit sphere
 $\{ x:\,\, |x|=1\}$
consists of  discrete  eigenvalues $\lambda_n$. Then
\begin{equation}\label{V3<LT}
\int_{|x|>1} |V(x)|^{1/2+p} |x|^{1-d}dx
\leq C\Bigl(1+\sum_n |\lambda_n|^{p}\Bigr),\qquad 0<p\leq 1/2,
\end{equation}
where the constant $C$ depends on $d$ and $p$ but does not  depend on $V$. Moreover,
\begin{equation}\label{V3<LT2}
\sum_{n=1}^\infty\Bigl(\int_{n<|x|\leq n+1} V(x) |x|^{1-d}dx\Bigr)^{2p}
\leq C_{\lambda_0}\Bigl(1+\sum_n |\lambda_n|^{p}\Bigr),\quad p\geq 1/2,
\end{equation}
where $C_{\lambda_0}$ depends on $d$,  $p$ and the lowest eigenvalue $\lambda_0$ of $H$.
\end{theorem} The proof of this theorem will be given in
Section 3, which  contains  a description of a simple trick reducing the case $d\geq2$  to the case $d=1$. In order  to complete the proof, it remains to  to note that the situation when $d=1$ was considered by Damanik and Remling  \cite{DR}.

Theorem~\ref{N} is related to the  case $d=2$. It is natural to ask what happens in higher dimensions? The following statement  aswers the question.

\bigskip

\begin{theorem}\label{smolchan} {\bf (Molchanov \cite{Mol}).}  Let $0<p<\infty$  and let $N$ be the number of negative    eigenvalues of $H=-\Delta-V$ acting in $L^2({\Bbb R}^d)$ with $d\geq3$. Then
there is no function   $F:\, {\Bbb R}_+\mapsto {\Bbb R}_+$,
 such that  
\begin{equation}\label{rlc}
\int_{{\Bbb R}^d} V^{p}(x)dx\leq F(N),
\end{equation}
for all positive 
$V\in C^\infty_0({\Bbb R}^d)$. In particular, one can find a smooth potential 
(with an infinite support), such that $N=0$ but $\int_{{\Bbb R}^d} |V|^p dx=\infty$ for all $p>0.$
\end{theorem}
\bigskip

Let  us discuss  the case   when $V$ changes its sign:
$$
V=V_+-V_-,\qquad 2V_\pm=|V|\pm V.
$$
Clearly, it is insufficient to consider only one Schr\"odinger operator
in this case. One has to treat $V$ and $-V$ symmetrically and  study  the spectra of two operators $H_+=-\Delta+V$ and $H_-=-\Delta-V$.
 \begin{theorem} \label{main} {\rm (see \cite{S})} Let
$V\in L^\infty({\Bbb R}^d)$
 be a  real  function.
Let  the  essential spectrum of both operators $-\Delta +V$ and
 $-\Delta - V$   be either positive or empty. Assume that
the negative eigenvalues  of the operators $-\Delta +V$ and
 $-\Delta - V$, denoted by $\lambda_n(V)$ and $\lambda_n(-V)$,
 satisfy  the condition
 $$
\sum_n \sqrt{|\lambda_n(V)|}+\sum_n\sqrt{|\lambda_n(-V)|}<\infty.
 $$
  Then
 $$
 V=V_0+{\rm div} ( A)+| A|^2
 $$
 where $V_0$ and    ${\rm div}A$  are locally bounded, $A$ is  continuous and
 has  locally square integrable  derivatives,
 $$
\int (|V_0|+| A|^2)|x|^{1-d}dx<\infty.
$$
\end{theorem}
Thus, if $\lambda_n(\pm V)$ are in $\ell^{1/2}$, then $V$  either oscillates or decays at the infinity. As a  consequence, we obtain the following  statement  about the  absolutely continuous spectrum of $H_\pm$.

\begin{theorem}\label{1.1} Let   $V\in L^\infty({\Bbb R}^d)$
 be a real  function. Assume that the negative spectra of the operators $H_+=-\Delta +V$ and
 $H_-=-\Delta - V$ consist only of
 eigenvalues,   denoted by $\lambda_n(V)$ and $\lambda_n(-V)$,
which satisfy  the condition
 $$
\sum_n \sqrt{|\lambda_n(V)|}+\sum_n\sqrt{|\lambda_n(-V)|}<\infty.
 $$
 Then  the   absolutely  continuous  spectra of  both operators
 are  essentially supported by $[0,+\infty)$, i.e.
the spectral projection $E_{H_\pm}(\delta)$ associated   to any subset 
$\delta\subset {\Bbb R}_+$ of
positive Lebesgue  measure is different  from zero.
\end{theorem}

\smallskip
While the proof  of this result is already   given in \cite{S}, we would  like to give it  again in order to
fill in the gaps and  make it more  clear.
Note  that this  theorem is proven in $d=1$ by Damanik and Remling \cite{DR}. One of the  missing parts  of the proof in  $d\geq 2$ was  the so called  trace formula obtained in 
\cite{LNS}. On the other hand, the application of  the main result of  \cite{LNS}  needs
a good  control of  negative  eigenvalues of the operators $H_\pm$. In particular, one has  to approximate  the potential $V$ by  
a compactly supported function in such a way
that  the corresponding eigenvalue sums $\sum_n |\lambda_n(\pm V)|^{1/2}$ would  not change much. It became possible  due to Theorem~\ref{main}, which gives an additional information about  $V$.

Below, we discuss  examples  of oscillating  potentials for which  the  corresponding eigenvalue sums are convergent.

\medskip

{\bf Example} (see \cite{Safronov}). Let $d\geq3$ and $V\in L^{d+1}({\Bbb R}^d)\cap
L^{\infty}({\Bbb R}^d)$ be a real potential whose Fourier
transform is  square integrable near the origin.  Then
$$
\sum_n \sqrt{|\lambda_n(V)|}+\sum_n\sqrt{|\lambda_n(-V)|}<\infty.
$$
In particular, the a.c. spectrum of $-\Delta+V$ is essentially  supported by the interval $[0,\infty)$.

\medskip

{\bf Example} (see \cite{SV}). Let $d\geq3$ 
and let $V_\omega$ be a real potential of the form
$$V_\omega(x)=\sum_{n\in {\Bbb Z}^d}\omega_n v_n\chi(x-n),$$ where $v_n$ are  fixed real    numbers, $\chi$ is the characteristic function of the  unit cube $[0,1)^d$
and $\omega_n$ are bounded independent   identically distributed random  variables with the zero expectation ${\Bbb E}(\omega_n)=0$.
If $V_\omega\in L^{d+2\gamma}({\Bbb R}^d)$ with $\gamma\geq 0$  for all $\omega=\{\omega_n\}$, then
$$
\sum_n |\lambda_n(V_\omega)|^\gamma<\infty,\qquad {\rm   almost \,\,\, surely.}
$$
In particular, if $V_\omega\in L^{d+1}({\Bbb R}^d)$  for all $\omega=\{\omega_n\}$, then the a.c. spectrum of $-\Delta+V_\omega$ is essentially  supported by the interval $[0,\infty)$.

\medskip

Note, that  one can not omit the  condition on the spectrum of one  of the operators $H_\pm$ in Theorem~\ref{1.1}, because 
the property of being  essentially supported by ${\Bbb R}_+$
 is not  shared by the spectra of all positive  Schr\"odinger operators.
One can conclude very little about the  spectrum  from the
fact  that  $-\Delta+ V\geq 0$, even if $V$ is bounded. For  instance, the  theory of
random  operators has  examples of   Schr\"odinger
operators with  positive $V$ and pure point spectra.
Therefore it is natural  to  combine the information given for
 $V$ and
$-V$. This idea was used in
\cite{DomK}  in dimension $d=1$ to prove the
following  striking result:
\begin{theorem}{\bf  (Damanik-Killip \cite{DomK})}\label{DomK} Let
$V\in L^\infty({\Bbb R})$.  If  the negative spectra  of the operators
$-\frac{d^2}{dx^2}+V$ and
 $-\frac{d^2}{dx^2}-V$ on the real line are finite,
 then  the positive spectra of these operators are  purely  absolutely  continuous.
\end{theorem}
It becomes clear from the nature  of the  obtained results that there is
a conservation law hiding behind  them. Appearance of negative
eigenvalues happens on the expense  of the absolutely continuous spectrum
which becomes "thinner". Therefore the thickness of the positive
spectrum can be estimated  by the number of  the negative energy levels. If
the number of these levels is "small", then the positive spectrum is
absolutely continuous. Put it differently, one can give a   quantitative formulation
of Theorem~\ref{1.1}. In order to  do that we have to introduce the spectral measure  of
the operator $H_+$ corresponding to an element $f$:
$$
\bigl((H_+-z)^{-1}f,f\bigr)=\int_{-\infty}^\infty \frac{d\mu(t)}{t-z},\qquad f\in L^2,
\quad \forall z \in {\Bbb C}\setminus {\Bbb R}.
$$
The integration in the right hands side is carried out with respect to  $\mu$ which is  the spectral measure of $H_+$.

\begin{theorem}\label{log} Let   $V\in L^\infty({\Bbb R}^d)$
 be a real  function. Assume that the negative spectra of the operators $H_+=-\Delta +V$ and
 $H_-=-\Delta - V$ consist only of
 eigenvalues,   denoted by $\lambda_n(V)$ and $\lambda_n(-V)$,
which satisfy  the condition
 $$
\sum_n \sqrt{|\lambda_n(V)|}+\sum_n\sqrt{|\lambda_n(-V)|}<\infty.
 $$
Then there exists an element $f\in L^2$ with $||f||=1$ such that  for any continuous compactly supported function $\phi\geq 0$ on the positive half-line $(0,\infty)$,
\begin{equation}\label{entr1}\begin{split}
\int_0^\infty\log \Bigl(\frac {\mu'(\lambda)}{\phi(\lambda)}\Bigr)\phi(\lambda)d\lambda\geq- C(\sum_n \sqrt{|\lambda_n(V)|}+\\
\sum_n\sqrt{|\lambda_n(-V)|}
+\sqrt{||V||_\infty}+1),
\end{split}
\end{equation}
where $||V||_\infty$ denotes the $L^\infty $-norm of the function $V$ and $C>0$  depends on $\phi$.
\end{theorem}

\bigskip

This theorem  implies  the  statement of Theorem~\ref{1.1}, because if the   eigenvalue sums  in the  right hand side of  \eqref{entr1} are convergent, then $\mu'(\lambda)>0$ almost  everywhere on the interval $[0,\infty)$.

Our  section "References" suggests a list of papers containing  the material  on both topics: eigenvalue estimates and  absolutely continuous spectrum of Schr\"odinger operators. 
All  papers \cite{1}-\cite{25} are  highly  recommended.

\section{Proof of  the estimate \ref{V<N2}}
 Our  arguments can be compared with the constructions
of the paper  \cite{DR} where similar questions  were studied for the case $d=1$.
The  following  statement is  obvious.
\begin{lemma}\label{NL1} Let $\Omega$ be a bounded domain in ${\Bbb R}^2$ with a piecewise smooth boundary.
 Let $\phi$ be a real valued bounded function with bounded  derivatives.
Suppose that $-\Delta \psi - V \psi=\lambda \psi$ and  the product $\phi\psi$ vanishes on the  boundary of
$\Omega$.
Then
$$ \int_{\Omega}\Bigl(|\nabla (\phi\psi)|^2-
V|\phi\psi|^2\Bigr)dx=\int_{\Omega}\Bigl(|\nabla \phi|^2\psi^2+\lambda |\phi\psi|^2\Bigr)dx.$$
\end{lemma}
Let  us introduce the unit square $Q=(-1,1)^2$ in the plane ${\Bbb R}^2$.

\begin{lemma}\label{NL2} Let $\Omega$ be an open domain with a piecewise smooth boundary.
Assume that the  lowest  eigenvalue $-\gamma^2$ of
$H$
on  $\Omega$ is  negative.
Then there  exists a
square $D$,
\begin{equation}\label{Omega}D=x_0+6 \gamma^{-1}Q,
\end{equation} such that
 $H$ restricted onto
$ \Omega\cap D$  has an eigenvalue not  bigger  than $-\gamma^2/2$
\end{lemma}

{\it Proof.} Let $\psi$ be  the  eigenfunction  corresponding  to the eigenvalue $-\gamma^2$
for the problem on the domain
$\Omega$ with the Dirichlet  boundary conditions. Put $L=2 \gamma^{-1}$ and pick a
point $x_0$  which gives  the maximum to the
functional $\int_{x_0+LQ}|\psi|^2dx$.  The latter integral is a continuous function of 
$x_0$,  
tending to zero as $|x_0|\to \infty$, so it does have a maximum. Suppose that  the coordinates of $x_0$ are the numbers $a$ and $b$, i.e. $x_0=(a,b)$. Define $$\phi(x)=\min\{\phi_0(\xi-a),\phi_0(\eta-b)\},\qquad x=(\xi,\eta),$$ where
\begin{equation}\label{Ndefphi}
\phi_0(t)=\begin{cases}
1,\quad {\rm if}\quad |t|<L,
\\0,\quad  {\rm if}\quad |t |\geq 3L,\\3/2-|t |/(2L), \qquad
{\rm otherwise.}
\end{cases}
\end{equation}
It is clear  that the support of $\phi$ is the square $D$ given by \eqref{Omega}.
 Now  the interesting  fact  is  that
$$
\int_{\Omega}|\nabla \phi|^2\psi^2dx\leq\frac{ \gamma^2}{ 2 }\int_{\Omega} |\phi\psi|^2dx
$$
which is  guaranteed  by the choice of $x_0$. Indeed, $|\nabla \phi|$ vanishes everywhere except for  eight  squares with the side length  $2L$, where  $|\nabla \phi|$  equals $1/(2 L) $.  Consequently,
$$
\int_{\Omega}|\nabla \phi|^2\psi^2dx\leq 2L^{-2}\int_{x_0+LQ} |\psi|^2dx.
$$
Therefore by Lemma~\ref{NL1}
$$ \int_{\Omega\cap D}\Bigl(|\nabla (\phi\psi)|^2+
V|\phi\psi|^2\Bigr)dx\leq-\frac{\gamma^2}{2}\int_{\Omega\cap D} |\phi\psi|^2dx.
$$
That proves  the statement. $\ \ \ \ \ \Box$
\bigskip

We  will also need  the following  result
\begin{proposition}\label{NA}
Let $H\geq -\gamma^2$ on a domain $\Omega$. Then
\begin{equation}\label{NL3}
\int_{\Omega}V(x)|\phi|^2dx\leq \Bigl(
\gamma^2\int_{\Omega} |\phi|^2dx+\int_{\Omega}|\nabla\phi|^2dx\Bigr)
\end{equation}
for any $\phi\in C_0^\infty(\Omega)$.  The  inequality \eqref{NL3} can be  extended to
  functions from the Sobolev  space $ W^{1,0}_2(\Omega)$.
\end{proposition}

The proof of the following  statement is a rather obvious consequence of Lemma ~\ref{NL2}

\begin{lemma}\label{NTL1} Suppose that the operator $H$ on the whole plane ${\Bbb R}^2$
has $N$ negative eigenvalues.
There is a  collection  of not more than $N$   squares $\Omega_n=x_n+L_nQ$ 
and not more than $N$ numbers $\epsilon_n>0$, such that
  $L_n=6\epsilon_n^{-1/2}$ and
\begin{equation}\label{geq1}
H\geq -\epsilon_n,\qquad {\rm on} \qquad {\Bbb R}^2\setminus\cup_{j<n}\Omega_j
\end{equation}
The collection of sets    fulfills
the condition, that $H\geq0$ on the set
${\Bbb R}^d\setminus\cup_n \Omega_n$. Moreover
$$
H\geq -\epsilon_{n}/2,\quad {\rm on}\,\, \Omega_n\setminus\cup_{j<n}\Omega_j.
$$
\end{lemma}

\bigskip

\bigskip

We  continue our  reasoning as follows.
Let $\phi_0$ be the function of  one real variable defined in \eqref{Ndefphi} with $L=1$. We introduce the functions $\psi_n$ on ${\Bbb R}^2$ by setting
$$
\psi_n(x)=\min\{\phi_0(L_n^{-1}(\xi-a_n)),\phi_0(L_n^{-1}(\eta-b_n))\},\qquad {\rm where}\,\,\,\, x=(\xi,\eta),\ \ \ x_n=(a_n,b_n).
$$
Note, that the  functions $\psi_n$  are equal to 1 on the cubes $\Omega_n$, but nevertheless  they are  compactly  supported.

Now we introduce another collection of functions  $\zeta_n$  equal to 1 on the union $\cup_{j\leq n}\Omega_j$.
Set  $\zeta_1(x)=\psi_1(x)$.
Given $\zeta_j$ for $j<n$, we  define the  functions $\zeta_n$ and $\tilde\psi_n$ as
$$
\zeta_n(x)=\max\{\zeta_{n-1},\psi_n(x)\},
$$
\begin{equation}\label{using}
\tilde\psi_n(x)=\min\{(1-\zeta _{n-1}(x)), \psi _n ((x-2x_n)/3)\}.
\end{equation}
It is easy to see that $\zeta_n=1$ on $\cup_{j\leq n}\Omega_j$ and    the support of
$\tilde\psi_n$  is contained in ${\Bbb R}^2\setminus\cup_{j<n}\Omega_j$.
Since the support of  $\tilde\psi_n$ is contained in ${\Bbb R}^2\setminus\cup_{j<n}\Omega_j$ and the operator $H$ on this set is larger  than $-\epsilon_n$,
we obtain that
$$
\int|\nabla\tilde\psi_n|^2dx-\int V(x)|\tilde\psi_n|^2dx\geq
-\epsilon_n\int |\tilde\psi_n|^2dx
$$
due to \eqref{geq1}.  Using \eqref{using} we derive the  estimate
\begin{equation}\label{intv}
\int V(x)|\tilde\psi_n|^2dx\leq \epsilon_n(6L_n)^2+\int|\nabla\zeta_{n-1}|^2dx+\int|\nabla\psi_n|^2dx,
\end{equation}
where the constant $C$ does not depend on $V$. Since   the  functions $\zeta_{n-1}$ and $\psi_n$ are  not larger  than 1,  it is easy to see that
$$
\int|\nabla\zeta_{n}|^2dx\leq\int|\nabla\zeta_{n-1}|^2dx+
\int|\nabla\psi_{n}|^2dx,
$$
which leads to the estimate
\begin{equation}\label{cn2}
\int|\nabla\zeta_{n}|^2dx\leq  8 n,\qquad \forall n\leq N.
\end{equation}
 Combining \eqref{intv} and \eqref{cn2}, we obtain that
$$
\int_{{\rm supp}\psi_n\setminus{\rm supp}\zeta_{n-1}} V(x) dx\leq\int V(x)|\tilde\psi_n|^2dx\leq 6^4+8n,
$$
On the other hand, it is  obvious that
$$
{\rm supp}\zeta_{n}={\rm supp}\psi_n\cup{\rm supp}\zeta_{n-1},\qquad \forall n\leq N,
$$
which leads to the estimate
\begin{equation}\label{intv2}
\int_{{\rm supp}\zeta_N} V(x) dx\leq\sum_n\int_{{\rm supp}\psi_n\setminus{\rm supp}\zeta_{n-1}} V(x) dx\leq (6^4+4(N+1))N,
\end{equation}

Define  now the function $\psi_{N+1}$ by
$$
\psi_{N+1}=\min\{\phi_0(\xi/R),\phi_0(\eta/R)\},\qquad x=(\xi,\eta)
$$
for $R$ large enough.  Then for $\tilde\psi=\min\{(1-\zeta_N),\psi_N\}$ we have the inequality
$$
\int|\nabla\tilde\psi|^2dx-\int V(x)|\tilde\psi|^2dx\geq
0,
$$
which implies that
\begin{equation}\label{final}
\int_{{\Bbb R}^d\setminus{\rm supp}\zeta_N} V(x) dx\leq\lim_{R\to\infty}\int|\nabla\tilde\psi|^2dx\leq 8(1+N).
\end{equation}
Actually, one can even get rid of 1 in the right hand side of \eqref{final} by considering a slightly different family  of  functions $\psi_{N+1}$.
Combining \eqref{intv2} with \eqref{final} we obtain \eqref{V<N2}.

\section{Proof of Theorem ~\ref{LTd3}}
Here we prove the estimates  \eqref{V3<LT} and \eqref{V3<LT2}. Let $N(\lambda, A)$ be the number  of eigenvalues of an operator $A=A^*$ situated to the left the point $\lambda\in {\Bbb R}$. Then according to the variational principle,
\begin{equation}\label{Var}\begin{split}
N(\lambda, A)=\max {\rm dim}  F\qquad \\ (Au,u)\leq \lambda||u||^2, \qquad  u\in F.
\end{split}
\end{equation}
It follows immediately from \eqref{Var}, that if $P_0$
denotes the orthogonal projection onto the space of spherically symmetric functions, then
$$
N(\lambda, P_0HP_0)\leq N(\lambda,H),\qquad \lambda\leq0,
$$
because,  we
 diminish  the number of considered subspaces by saying that $F$ in \eqref{Var} consists only of  spherically symmetric functions.
 
 On the other side,  $P_0HP_0$ is unitary equivalent to the orthogonal sum of the zero operator and the one-dimensional Schr\"odinger operator
 $$
 -\frac{d^2}{dr^2}+\frac{\alpha_d}{ r^2} -\bar V(r),
 $$
where 
 $$
 \bar V(r)=\frac1{c_d r^{d-1}}\int_{|x|=r} V(x)
dS $$
and $c_d$ is the  area of the unit sphere. It remains to note that Theorem~\eqref{LTd3} is proven for $d=1$ in \cite{DR}.

\section{Proof of Theorem ~\ref{main}}

Following the main idea of  \cite{DR}, we find the regions where the  eigenfunctions of $H_\pm$
live. Such regions have  the property that after introducing the Dirichlet condition on their boundary, the corresponding negative eigenvalue  does not get  too close to zero.
The  next statement is  obvious.
\begin{lemma}\label{3L1}
Let $\phi$ be a real valued bounded function with bounded  derivatives.
Suppose that $-\Delta \psi \pm V \psi=\lambda \psi$ and  the product $\phi\psi$ vanishes on the  boundary of  the  domain
$\{a<|x|<b\}$.
Then
$$ \int_{a<|x|<b}\Bigl(|\nabla (\phi\psi)|^2\pm
V|\phi\psi|^2\Bigr)dx=\int_{a<|x|<b}\Bigl(|\nabla \phi|^2\psi^2+\lambda |\phi\psi|^2\Bigr)dx$$
\end{lemma}
Before stating a very important  lemma we introduce the  notion of the internal size (width) of a
 spherical layer $\{a<|x|<b\}$, which is  $b-a$. We will also say that one set $X\subset {\Bbb R}^d$ is situated to the "left" of another $Y\subset {\Bbb R}^d$ (correspondingly, $Y$ is  situated to the "right" of $X$) if $|x|<|y|$ for all $x\in X$ and $y\in Y$.
\begin{lemma}\label{3L2} Assume that the  lowest  eigenvalue $-\gamma^2$ of $H_\pm$
on the domain $\{a<|x|<b\}$ is  negative.
If $b-a\geq 6\gamma^{-1}$, then there  exist a
spherical layer $\Omega\subset\{a<|x|<b\}$ whose  width is $d(\Omega)=6\gamma^{-1}$
such that $H_\pm$ restricted onto
$ \Omega$  has an eigenvalue not  bigger  than $-\gamma^2/2$.
\end{lemma}

{\it Proof.} Let $\psi$ be  the  eigenfunction  corresponding  to the eigenvalue $-\gamma^2$
for the problem on the domain
$\{a<|x|<b\}$ with the Dirichlet  boundary conditions. Put $L=\gamma^{-1}$ and pick a
number $c>0$  which gives  the maximum to the
functional $\int_{c-L<|x|<c+L}|\psi|^2dx$.  The latter functional is a continuous function of $c$,  tending to zero as $c\to \infty$, so it does have a maximum. Define
\begin{equation}\label{defphi}
\phi(x)=\begin{cases}
1,\quad {\rm if}\quad ||x|-c|<L, \\0,\quad  {\rm if}\quad ||x|-c|\geq 3L,\\3/2-||x|-c|/(2L), \qquad
{\rm otherwise.}
\end{cases}
\end{equation}
Let $\Omega$  be the intersection of the support of $\phi$ with $\{a<|x|<b\}$.  Without
loss  of generality we  can assume that
$d(\Omega)=6\gamma^{-1}$. Otherwise, we can make it larger, so that the bottom of the spectrum will not increase.
 Now  the interesting  fact  is  that
$$
\int_{a<|x|<b}|\nabla \phi|^2\psi^2dx\leq\frac{ \gamma^2}{ 2 }\int_{a<|x|<b} |\phi\psi|^2dx
$$
which is  guaranteed  by the choice of $c$. Indeed, $|\nabla \phi|$ vanishes everywhere except for  two  spherical layers of width $2L$, where  it equals $\gamma/2 $.  Consequently,
$$
\int_{a<|x|<b}|\nabla \phi|^2\psi^2dx\leq \frac{ \gamma^2}{ 2 }\int_{c-L<|x|<c+L} |\psi|^2dx.
$$
Therefore by Lemma~\ref{3L1}
$$ \int_{a<|x|<b}\Bigl(|\nabla (\phi\psi)|^2\pm
V|\phi\psi|^2\Bigr)dx\leq-\frac{\gamma^2}{2}\int_{a<|x|<b} |\phi\psi|^2dx.
$$
That proves  the result. $\ \ \ \ \ \Box$

\bigskip

We  also need  the following elementary statement:
\begin{lemma}\label{A}
Let $H_\pm\geq -\gamma^2$  on a bounded spherical  layer $\Omega=\{a<|x|<b\}$, $a>0$,   for  both  indexes $\pm$. Then
$V +\gamma^2={\rm div} A+|A|^2$ on $\Omega$,
where $A$  satisfies  the estimate
\begin{equation}\label{L3}
\int_{a<|x|<b}|\phi|^2|A(x)|^2dx\leq C\Bigl(
\gamma^2\int_{a<|x|<b} |\phi|^2dx+\int_{a<|x|<b}|\nabla\phi|^2dx\Bigr)
\end{equation}
for any $\phi\in C_0^\infty(\Omega)$ with a  constant $C$ independent of $\gamma$, $\Omega$ and $\phi$.
\end{lemma}

{\it Proof.} The representation $V =-\gamma^2+{\rm div} A+|A|^2$ on $\Omega$ follows  from the results
of \cite{DKSaf}. The idea is   to define the vector potential by
$
A=u^{-1}\nabla u,$
  where $u$ is a positive solution of the equation  $ -\Delta u+Vu =\gamma^2 u.$
Now
$$
\int_{a<|x|<b}\Bigl(|\nabla \phi|^2
-V|\phi|^2\Bigr)dx\geq -\gamma^2\int_{a<|x|<b}|\phi|^2dx
$$
which leads  to the inequality
\eqref{L3} due to the estimate
$$
\int_{a<|x|<b} {\rm div} A |\phi|^2dx\leq \epsilon \int_{a<|x|<b} |A|^2|\phi|^2 dx+
\epsilon^{-1}\int_{a<|x|<b}|\nabla\phi|^2dx.
$$
The proof is completed. $\ \ \ \ \ \ \Box$
\smallskip

The main ingredients of our proof  are in the following technical lemmas, which
can be compared  with the corresponding  set of  statements  from \cite{DR}. Our  proof  is  shorter
because  instead of functions   with symmetric  graphs  we  will  use  functions whose  graphs will have three  different  
slopes.
This will influence the choice of  sets $\Omega_n$.
\begin{lemma}\label{TL1} Let $V(x)=0$ for $|x|<2$.
There is a  sequence of spherical layers $\Omega_n$
and a sequence of numbers $\epsilon_n>0$, such that $\sum_n
\epsilon_n^{1/2}<\infty$
and the width of  $\Omega_n$
is bounded by $C\epsilon_n^{-1/2}$ with some $C$ independent of $n$. The sequence of sets    fulfills
the condition, that $H_\pm\geq0$ on the set
${\Bbb R}^d\setminus\cup_n \Omega_n$. Moreover
$$
H_\pm\geq -\epsilon_{j(n)},\quad {\rm on}\,\, \Omega_n
$$
where $j(n)$ is the lowest number $j$ such that $\Omega_j\cap\Omega_n\neq\emptyset$. If
$\Omega_{j}\cap\Omega_n\neq\emptyset$
and the width of $\Omega_{j}\cap\Omega_n$ is  bounded  from below by $6(1-20^{-1})\epsilon_{k}^{-1/2}$:
$$
d\Bigl(\Omega_{j}\cap\Omega_n\Bigr)\geq 6(1-20^{-1})\epsilon_{k}^{-1/2},
$$
 where $k={\rm min}\{j,n\}$. The  choice of the sequences can be so made that
 for each $n$  the number of the sets $\Omega_j$  intersecting 
 $\Omega_n$ is not bigger than 2 and 
 $$
{\rm dist}\Bigl(\Omega_n,\cup_{m<j(n)}\Omega_m\Bigr)\geq\frac{3}{10\epsilon^{1/2}_{j(n)}}.
$$
\end{lemma}
\bigskip
{\it Proof.}
In the proof,  we  also need  to construct some sequence of sets $\omega_n\subset \Omega_n$. Put $\Omega_0=B_2$
and $\omega_0=B_1$ where $B_r$ denotes  the ball of radius $r>0$  about  the  origin; $\epsilon_0$
can be any  sufficiently large  number, having the property that $-\epsilon_0$
lies below the bottom of the spectra of operators $H_\pm$.

Given  $\omega_n\subset\Omega_n$ and $\epsilon_n$ for $n<N$ we  consider the set
$$
S={\Bbb R}^d\setminus\cup_{n<N} \Omega_n
$$
and  define $-\epsilon_N$ as  the  lowest eigenvalue
of both operators $H_\pm$ on $S$. Note that
it follows    immediately that
$$
\epsilon_j\geq\epsilon_{j+1}.
$$
Let
$\omega_N\subset S$ be the largest spherical layer
where one of the operators $H_\pm$ has spectrum below $-\epsilon_N/2$, i.e.
$$
\inf \sigma(H_\pm)\leq -\epsilon_N/2,
$$
and the width of $\omega_N$ is not bigger than $L= 6 \epsilon_N^{-1/2}$. The  existence of this set is proven in Lemma 2.2. We  would like to  stress that there  are  two options:

\bigskip
\noindent

1) either the boundary of $\omega_N$ contains at least one  interior point  of  $S$  or

2)  the boundary of $\omega_N$ is contained in   the boundary   of  $S$.

 \bigskip
\noindent
In the first case,  the width of $\omega_N$ is equal to $L= 6 \epsilon_N^{-1/2}$.

Denote by $S_+$ and  $S_-$ the connected component of
$S\setminus \omega_N$ situated to the right and  to the left of $\omega_N$
correspondingly. In the second of the  above cases, both  sets $S_\pm$ are  empty.
Let $\Omega_j=:\Omega_-$
and $\Omega_k=:\Omega_+$, $j,k<N$,
 be the  two sets  which have common boundary with $S_-$ and $S_+$
correspondingly. Denote  $\omega_-=\omega_j$
and  $\omega_+=\omega_k$. Our  construction (or induction assumptions)
allow  us to assume that $${\rm dist}\{\omega_\pm,S_\pm\}\geq L_\pm,\,{\rm where}\, L_-=6 \epsilon_j^{-1/2}\,{\rm and}\,
L_+=6 \epsilon_k^{-1/2}$$
If the width of $S_\pm$  is not bigger than $3L$ we include
$\Omega_\pm\setminus\{x:\,\,{\rm dist}(x,\omega_\pm)\leq L_\pm/20\}$  into $\Omega_N$
by definition by
 demanding that $\{x:\,\,{\rm dist}(x,\omega_\pm)\leq L_\pm/20\}$ and $\Omega_N$ has
a non- empty piece of common boundary and $\{x:\,\,{\rm dist}(x,\omega_\pm)\leq L_\pm/20\}\cap\Omega_N=\emptyset$.
Otherwise,
$$S_\pm\setminus \{x: \ {\rm dist}(x,\omega_N)<L\}=S_\pm\setminus\Omega_N.$$
Obviously $\Omega_N$ is contained in ${\Bbb R}^d\setminus\cup_{i< j(N)} \Omega_i$
where $j(N)$ is the lowest number $i$ such that $\Omega_i\cap\Omega_N\neq\emptyset$.
Therefore,
$$
H_\pm\geq- \epsilon_{j(N)}\quad{\rm on}\quad \Omega_N.
$$

Observe that the distance
from the boundary of $\Omega_N$ to $\omega_N$ is not less  than the smallest of the numbers $19L_\pm/20$, which equals $57/(10\sqrt{\epsilon_{j(N)}})$.
Obviously, for any $\gamma>0$ there exist  a number $N$ such that the infinum of the spectrum of
both operators
$H_{\pm}$ on the domain
$$
{\Bbb R}^d\setminus\cup_{n<N} \Omega_n
$$  is  higher  than  $-\gamma$.
Assume the  opposite. Then for any $N$ one of the   operators $H_{\pm}$ on the domain
$$
{\Bbb R}^d\setminus\cup_{n<N} \Omega_n
$$  has an eigenvalue which is not bigger   than  $-\gamma$. Then
there is an eigenvalue of  one of the the operators on  $\omega_N$ which is  not  bigger than
$-\gamma/2$.
This  implies  that the negative spectrum of one of  the operators $H_\pm$ is not  discrete. So we come to the  conclusion that $H_\pm>0$
on
$$
{\Bbb R}^d\setminus\cup_{n} \Omega_n.
$$
Now  let  us observe  that $\sum_n \epsilon^{1/2}_n<\infty$, because the domains $\omega_n$ are disjoint and  our operators have an eigenvalue  lower than $-\epsilon_N/2$  on each of  $\omega_N$.
Also, it is  clear that any bounded ball $B_r$ of radius $r<\infty$ intersects only finite number of
$\Omega_n$, otherwise  a Schr\"odinger operator on $B_r$ would
have infinite number of  eigenvalues below zero.
$\,\,\,\,\Box$

\bigskip

Lemma~\ref{TL1} allows one to estimate  the potential on  the union of the  sets $\Omega_n$
However, these sets might not exhaust the whole space  ${\Bbb R}^d$, so we have to consider the case when
$$
{\Bbb R}^d\setminus\cup_{n} \Omega_n\neq\emptyset.
$$
The set ${\Bbb R}^d\setminus\cup_{n} \Omega_n$ might contain unfilled gaps, where   both operators $H_+$ and $H_-$ are positive. Suppose that a spherical layer $\Lambda$ is a connected  component of the set  ${\Bbb R}^d\setminus\cup_{n} \Omega_n$. Then either  the width of 
the layer  $\Lambda$ is too large or it is small relative to the widths of  the sets  $\Omega_{n_1}$ and $\Omega_{n_2}$  having a common  boundary  with  the set $\Lambda$.
So we can  formulate
\begin{lemma}\label{TL2} Enlarging some of the sets $\Omega_n$ in
  Lemma~\ref{TL1} one can   achieve that 
\begin{equation} \label{z1}
{\Bbb R}^d=\Bigl(\cup_{n} \Omega_n\Bigr)\cup \Bigl(\cup_{n} \Lambda_n\Bigr).
\end{equation}
 where  $\Lambda_n$ are such spherical layers on which  both operators $H_+$ and $H_-$ are positive. 
Each $\Lambda_m$ intersects exactly  two  sets $\Omega_n$.
If $\Lambda_n$ intersects  $\Omega_{n_1}$ and $\Omega_{n_2}$, then 
$$
d(\Lambda_n\cap \Omega_{n_j})\geq 5\epsilon^{-1/2}_{n_j},\qquad j=1,2,
$$
where $d(G)$ denotes the width of $G$. Moreover, $\Lambda_n$ in \eqref{z1} can be so chosen that
both operators $H_+$ and $H_-$ are positive on
$$
\Lambda_n\cup \Bigl( \cup_{j=1,2}\{x\in \Omega_{n_j}:\quad {\rm dist}(x,\Lambda_n)<\epsilon^{-1/2}_{n_j}\}\Bigr).
$$
\end{lemma}
{\it Proof}. Indeed,  each time  when  there is a  spherical layer $\Lambda$
 having a  common boundary with   $\Omega_{n_1}$ and $\Omega_{n_2}$ and  the property that $H_\pm\geq0$
on $\Lambda$, we  compare $d(\Lambda)$ with  the  numbers $6\epsilon^{-1/2}_{n_1}$ and $6\epsilon^{-1/2}_{n_2}$. If  $d(\Lambda)$ is  small (i.e. 
 $d(\Lambda)\leq 6(\epsilon^{-1/2}_{n_1}+\epsilon^{-1/2}_{n_2})$) we keep enlarging 
both sets $\Omega_{n_1}$ and $\Omega_{n_2}$ until $\Lambda$ disappears.
However if $d(\Lambda)> 6(\epsilon^{-1/2}_{n_1}+\epsilon^{-1/2}_{n_2})$, then we enlarge
$\Omega_{n_1}$ and $\Omega_{n_2}$ giving them the pieces of $\Lambda$ of the width 
$6\epsilon^{-1/2}_{n_1}$ and $6\epsilon^{-1/2}_{n_2}$, correspondingly.
Finally, we make the set $\Lambda$ smaller by removing two pieces of the width $\epsilon^{-1/2}_{n_1}$ and 
$\epsilon^{-1/2}_{n_2}$.
 $\ \ \ \ \ \Box$

 \begin{lemma} Let  $\Omega_n$,  $\epsilon_n$ and $\Lambda_n$ be  the same as in Lemma~\ref{TL1} and 
Lemma~\ref{TL2}. There  exists a   sequence of $H^1$-functions  $\phi_n\geq0$
supported in  $\Omega_n$ and  a   sequence of $H^1$-functions  $\psi_n\geq0$
supported in  $\Lambda_n$ such  that
\begin{equation*}\begin{split}
\sum_n\phi_n(x)+\sum_n\psi_n(x)=1,\,\,\,\,\,\,\,\,\,\,\,\,\,\,\,\,\,\,\,\,\,\,\,\,\,\,\, \\ \sum_n\int (|\nabla\phi_n(x)|^2+  |\nabla\psi_n(x)|^2)|x|^{1-d}dx
\leq C\sum_n
\epsilon_n^{1/2}.
\end{split}
\end{equation*}
Moreover, one can find vector potentials $A_n$ and $\tilde A_n$ such that 
\begin{equation}\label{reprV}
V+\epsilon_{j(n)}={\rm div} A_n+|A_n|^2\,\, {\rm on}\,\,  \Omega_n,\qquad  
V={\rm div} \tilde A_n+|\tilde A_n|^2\,\, {\rm on} \,\, \Lambda_n
\end{equation}
and
$$\sum_n\Bigl(\int_{\Omega_n}|A_n|^2|x|^{1-d}dx+ \int_{\Lambda_n}|\tilde A_n|^2|x|^{1-d}dx\Bigr)\leq C\Bigl(1+\sum_n
\epsilon_n^{1/2}\Bigr).$$
\end{lemma}
\bigskip
{\it Proof.} According to our constructions,  if $\Omega_{j}\cap\Omega_n\neq \emptyset$  
then the width of $\Omega_{j}\cap\Omega_n$ is  bounded  from below by $6(1-20^{-1})\epsilon_{k}^{-1/2}$, where $k={\rm min}\{j,n\}$, while 
 the width of the intersections 
$
\Lambda_i\cap \Omega_{n}\neq \emptyset$ 
are  not less than 
$ 
5\epsilon^{-1/2}_{n}.
$
We define the functions $\phi_n$ in such a way that the slopes of their graphs  will be  inversely  proportional to the width
of these intersections.
Therefore 
 the functions $\phi_n$ will satisfy the  relations
$$
\int_{\Omega_j\cap \Omega_n}|\nabla\phi_n|^2|x|^{1-d}dx\leq C \epsilon_k^{1/2},
\qquad k={\rm min}\{j,n\};
$$
$$
\int_{\Omega_n\cup\Lambda_j}|\nabla\phi_n|^2|x|^{1-d}dx\leq C \epsilon_n^{1/2}.
$$
 We can additionally require  that  $\phi_n+\phi_j=1$ on the intersection $\Omega_n\cap\Omega_j$ and
$\phi_n+\psi_i=1$ on the intersection $\Omega_n\cap\Lambda_i$. 
The representations \eqref{reprV} could be obtained by application of Lemma~\ref{A}, since the operators $H_\pm$ satisfy the inequalities
$$
H_\pm\geq -\epsilon_{j(n)}\,\, {\rm on}\,\,  {\Bbb R}^d\setminus \Bigl(\cup_{m<j(n)}\Omega_m\Bigr),$$ and since both operators
$H_+$ and $H_-$ are positive on
$$
\Lambda_n\cup \Bigl( \cup_{j=1,2}\{x\in \Omega_{n_j}:\quad {\rm dist}(x,\Lambda_n)<\epsilon^{-1/2}_{n_j}\}\Bigr),
$$ where  $\Omega_{n_j}\cap \Lambda_n\neq \emptyset.$
The estimates for $A_n$ in this  construction  follow  from Lemma~\ref{A} with $\phi=|x|^{(1-d)/2}\tilde{\phi_n}$, where
$\tilde{\phi_n}=1$ on the set $$\Omega_n\subset {\Bbb R}^d\setminus \Bigl(\cup_{m<j(n)}\Omega_m\Bigr).$$ The slopes of the remaining part of the graph of this function should be inversely  proportional to the number
$$
{\rm dist}\Bigl(\Omega_n,\cup_{m<j(n)}\Omega_m\Bigr)\geq\frac{3}{10\epsilon^{1/2}_{j(n)}}.
$$
The latter inequality was mentioned in  Lemma~\ref{TL1}.
The integral estimates for the vector potentials $\tilde A_n$ can be obtained  in a similar way. $\ \ \ \ \ \Box$

\bigskip
{\it  The end  of the proof of Theorem~\ref{main}.} Let us define
$$
A=\sum_n (\phi_n A_n+\psi_n \tilde A_n),\qquad W=-\sum_n \epsilon_{j(n)}\phi_n
\quad
 V_1=
 W+{\rm div} ( A)+| A|^2
$$
Then one can easily
see  that
$$
V_1=
V+\sum_n ( A_n\nabla\phi_n +\tilde A_n\nabla\psi_n) +|A|^2 -\sum_n 
(\phi_n| A_n|^2+ \psi_n| \tilde A_n|^2),$$
 which   implies 
  $$
\int |V-V_1||x|^{1-d}dx<\infty.
$$
If we define $V_0=V-V_1+W$, then we will represent  the potential in the form
 $V=V_0+{\rm div} ( A)+| A|^2$. 
Finally, since
$$
\int |W||x|^{1-d}dx<\infty,
\quad \int |A|^2|x|^{1-d}dx<\infty,
$$
 we conclude that
$$
\int (|V_0|+| A|^2)|x|^{1-d}dx<\infty.\,\,\,\,\,\,\Box
$$

\bigskip

\section{ Proof of Theorem ~\ref{1.1}}
Now, we use the results of \cite{LNS} to  establish the  presence of the   absolutely continuous spectrum of $H_+$. As it is known the operator $-\Delta$ is unitary equivalent to
the sum $-\partial^2/\partial r^2+\alpha_d/r^2+(-\Delta_\theta)/r^2 $ where $\Delta_\theta$
is the Laplace-Beltrami operator on the unit sphere and $\alpha_d$ is a certain constant.
Since the values of $V$ on compact subsets do not influence the presence of absolutely continuous spectrum, 
we can assume that $V(x)=-\alpha_d/|x|^2$ for $1<|x|<2$. According to \cite{LNS},
there  is  a probability measure $\mu$ on the real line
${\Bbb R}$, whose a.c. component is  essentially supported  by a  subset  of  the a.c. spectrum of the
operator $H_+$. Namely, one constructs an operator $A_+$  having the same a.c. spectrum as $H_+$   in the following way:
$$
A_+=-\Delta+V,\qquad D(A_+)=\{u\in H^2({\Bbb R}^d \setminus B_1):\ \ \ \ u(\theta)=0,\ \  \theta
\in {\Bbb S}^{d-1}\}.
$$
Then one sets $\mu(\delta)=(E_{A_+}(\delta)f,f)$
for a spherically symmetric  function $f$ having the property
$$
{\rm supp}f\subset\{ x\in {\Bbb R}^d:  \quad 1<|x|<2\}.
$$
\begin{theorem}
{\bf \cite{LNS}}
Let  $V$ be a compactly  supported potential and let $\lambda_j(V)$ be the negative eigenvalues of $H_+$. Then for any continuous compactly supported function $\phi\geq 0$ on the positive half-line $(0,\infty)$, we have
\begin{equation}\label{entr}\begin{split}
\int_0^\infty\log \Bigl(\frac {\mu'(\lambda)}{\phi(\lambda)}\Bigr)\phi(\lambda)d\lambda\geq- C(\sum_j\sqrt{|\lambda_j(V)|}\\
+\sqrt{||V||_\infty}+\int V(x)|x|^{1-d}dx+1)
\end{split}
\end{equation}
where $||V||_\infty$ denotes the $L^\infty $-norm of the function $V$ and $C>0$  depends on $\phi$ and the choice of $f\in L^2$.
\end{theorem}

One of the important  properties of the measure $\mu$
is that (see \cite{Stolz})
$$
V_n\to V \ {\rm in}\ L^2_{loc}\Rightarrow \mu_n\to \mu \ \ {\rm weakly}.
$$ 
Using  are going to combine this property with the  lower semi-continuity  of the entropy (see the paper \cite{KS} by Killip and Simon)
\begin{theorem}\label{semicont }
Let the sequence of probability measures $\mu_n$ converge to $\mu$ weakly on the  real line ${\Bbb R}$.
Then for any $0<a<b<\infty$ and any  positive  continuous  (on the real line) function $\phi$ with the property ${\rm supp}\phi\subset [a,b]$, we have
$$
\int_a^b\log \Bigl(\frac {\mu'(\lambda)}{\phi(\lambda)}\Bigr)\phi(\lambda)d\lambda\geq \liminf_{n\to \infty} 
\int_a^b\log \Bigl(\frac {\mu'_n(\lambda)}{\phi(\lambda)}\Bigr)\phi(\lambda)d\lambda. $$
\end{theorem}
Let $[V_0]_+$ and  $[V_0]_-$   be the positive and negative   parts of the function $V_0$ and  let $\chi_n$ be  the characteristic  function of
the  ball of radius $n$ with the center at the point $0$.
Consider the Schr\"odinger operator with the potential $$V_{n,\varepsilon}=[V_0]_+-\chi_n[V_0]_-+{\rm div} ( A)+(1-\varepsilon)^{-1}| A|^2,\qquad \varepsilon\in(0,1).$$ 
The number of negative eigenvalues of the Schr\"odinger operator with the potential $V_{n,\varepsilon}$
 can be estimated by a   quantity that depends only on $\varepsilon^{-1}\chi_n[V_0]_-$.
Indeed,  define $N(W)$ as the number of negative eigenvalues of the operator $-\Delta+W$, where $W$ is a real potential.
Then 
$$
N(W_1+W_2)\leq N(\varepsilon^{-1}W_1)+N((1-\varepsilon)^{-1}W_2).
$$
Set now $W_1= -\chi_n[V_0]_-$
and $W_2=[V_0]_++{\rm div} ( A)+(1-\varepsilon)^{-1}| A|^2$.
Since $N((1-\varepsilon)^{-1}W_2)=0$, we will obtain that
$$
N(V_{n,\varepsilon})\leq N(-\varepsilon^{-1}\chi_n[V_0]_-).
$$
The right hand side of this inequality is independent of $A$ and $[V_0]_+$.
Therefore both of these functions can be  approximated by compactly  supported   functions in such a way that $N(V_{n,\varepsilon})$ stays bounded.
Let $\mu_{n,\varepsilon}$
be the corresponding spectral  measure of 
 the  operator  $-\Delta+V_{n,\varepsilon}$.  Then by the  lower semi-continuity of
the  entropy (see \cite{KS})
we obtain
\begin{equation}\label{semi}\begin{split}
\int_0^\infty\log \Bigl(\frac {\mu'_{n,\varepsilon}(\lambda)}{\phi(\lambda)}\Bigr)\phi(\lambda)d\lambda\geq- C(\sum _j|\lambda_j(V_{n,\varepsilon})|^{1/2}+\sqrt{||V||_\infty}\\ +
\int (|V_0|+(1-\varepsilon)^{-1}| A|^2)|x|^{1-d}dx+1)
\end{split}
\end{equation}
for any continuous   function $\phi\geq 0$ with supp$\phi\subset (0,\infty)$.
Since $V_{n,\varepsilon}$ is monotone in  both $n$ and $\varepsilon$, we have
$\lambda_j(V) \leq \lambda_j(V_{n,\varepsilon})$. 
Therefore,
the relation 
\eqref{semi} leads to the inequality
\begin{equation*}\begin{split}
\int_0^\infty\log \Bigl(\frac {\mu'_{n,\varepsilon}(\lambda)}{\phi(\lambda)}\Bigr)\phi(\lambda)d\lambda\geq-C(\sum _j|\lambda_j(V)|^{1/2}+\sqrt{||V||_\infty}\\ +
\int (|V_0|+(1-\varepsilon)^{-1}| A|^2)|x|^{1-d}dx+1).
\end{split}
\end{equation*}

Now  we can extend \eqref{entr} to  the general case due to the upper semi-continuity of the entropy:
\begin{equation*}\begin{split}\int_0^\infty\log \Bigl(\frac {\mu'(\lambda)}{\phi(\lambda)}\Bigr)\phi(\lambda)d\lambda\geq-
C(\sum_j\sqrt{|\lambda_j(V)|}+\sqrt{||V||_\infty}\\+\int (|V_0|+| A|^2)|x|^{1-d}dx+1).\end{split}
\end{equation*}
Convergence of the integral in the left hand side implies that $\mu'>0$ almost everywhere on the set where  $\phi$ is positive.
This completes  the proof  of both  Theorems~\ref{1.1} and ~\ref{log}.$\ \ \ \ \Box $

\end{document}